\documentclass[12pt]{article}

\newcommand{\ignore}[1]{}

\usepackage{amssymb}
\usepackage{amsthm}
\usepackage{graphicx}

\newcommand{\enum}[1]{\begin{enumerate} {#1} \end{enumerate}}

\newcommand{\dfn}[1]{\textbf{#1}}

\newcommand{\rr}{{\mathbb{R}}}

\newcommand{\ty}{\nabla\mathrm{Y}}
\newcommand{\yt}{\mathrm{Y}\nabla}

\newtheorem{theorem}{Theorem}
\newcommand{\thm}[1]{ \begin{theorem} {#1} \end{theorem} }
\newtheorem{corollary}[theorem]{Corollary}

\newtheorem{lemma}[theorem]{Lemma}

\theoremstyle{definition} 
\newtheorem{definition}{Definition}

\newtheorem{example}{Example}

\theoremstyle{plain} 

\newcommand{\pf}[1]{\begin{proof}{#1}\end{proof}}

\newtheorem*{theoremx}{Theorem}

\newtheorem*{lemmax}{Lemma}

\newtheorem*{definitionx}{Definition}

\newcounter{arb}

\begin{document}{

\title{The Y-triangle move does not preserve\\
Intrinsic Knottedness}
\author{Erica Flapan and Ramin Naimi}
\date{}

\maketitle

\abstract{
We answer the question
``Does the Y-triangle move preserve intrinsic knottedness?"
in the negative
by giving an example of a graph
that is obtained from the intrinsically knotted graph $K_7$
by triangle-Y and Y-triangle moves
but is not intrinsically knotted.
}

\section{Introduction}

A graph is said to be \dfn{intrinsically knotted}
(IK)
if every embedding of it in $\rr^3$
contains a cycle that is a nontrivial knot.
Similarly, a graph is said to be \dfn{intrinsically linked}
(IL)
if every embedding of it in $\rr^3$
contains a nontrivial link.
Sachs~\cite{Sa}
and
Conway and Gordon~\cite{CG}
showed that $K_6$,
the complete graph on six vertices,
is IL.
Conway and Gordon~\cite{CG} also
showed that $K_7$ is IK.

A \dfn{$\ty$ move} on an abstract graph
consists of removing the edges of a 3-cycle $abc$ in the graph,
and then adding a new vertex $v$ and connecting it
to each of the vertices $a$, $b$, and $c$,
as shown in Figure~\ref{triangleYmoveFigure}.
The reverse of this move is called a
\dfn{$\yt$ move}.
Note that in a $\yt$ move,
the vertex $v$ cannot have degree greater than three.

\begin{figure}[ht]

 \centering
 \includegraphics[width=80mm]{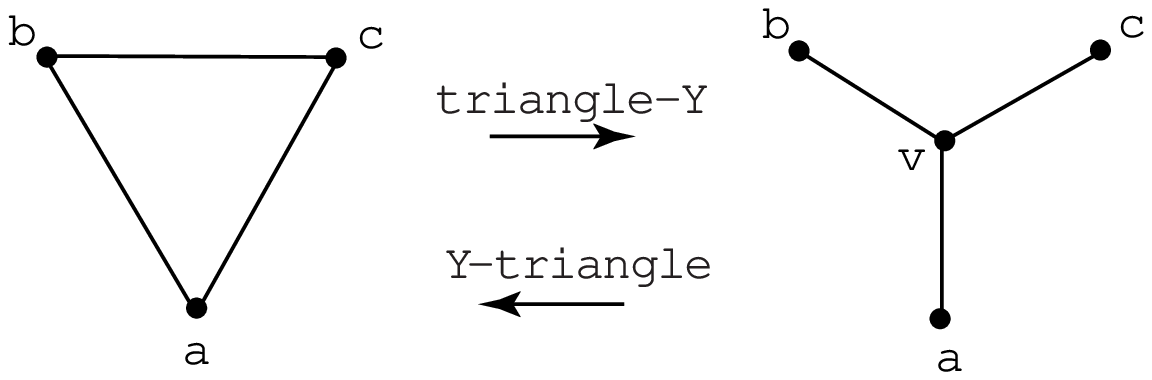}
 \caption{$\ty$ and $\yt$ moves.
 \label{triangleYmoveFigure}
 }
\end{figure}

Sachs~\cite{Sa} noticed that
additional IL graphs can be obtained from $K_6$
by doing finite sequences of $\ty$ and $\yt$ moves on it.
Motwani, Raghunathan, and Saran~\cite{MRS}
showed that performing a $\ty$ move
on any IK or IL graph
produces a graph with the same property.
Robertson, Seymour, and Thomas~\cite{RST} (Lemmas~1.2 and 5.1(iii))
proved that a $\yt$ move
on any IL graph produces an IL graph again.

It has been an open question
whether a $\yt$ move on an IK graph
always produces an IK graph again.
We prove that the answer is negative,
by giving a knotless embedding of a graph $G_7$
that is obtained from $K_7$
by $\ty$ and $\yt$ moves.

A graph $H$ is a \dfn{minor} of another graph $G$
if $H$ can be obtained from $G$
by a finite sequence of edge deletions and contractions
\cite{diestel}.
A graph is said to be \dfn{minor minimal} with respect to a property
if the graph has that property but
no minor of it has the property.

We work with connected, finite, simple graphs,
i.e., graphs with no loops (an edge whose endpoints are the same)
and no double-edges (two edges with the same pair of endpoints).
This is because loops and double-edges
do not affect whether or not a graph is IK or IL:
they can always be embedded such that they bound small disks
with interiors disjoint from the rest of the graph.
Thus, in edge contractions and $\yt$ moves on an abstract graph,
whenever a double-edge is introduced,
one of the two edges is deleted.

\section{Description of the graph $G_7$}

We label the seven vertices of the abstract graph $K_7$
with the letters $a$ through $g$.
We perform the following
five $\ty$ and two $\yt$ moves on $G_0 = K_7$
to obtain the graph $G_7$.

\enum{
\item
$G_0 \to G_1$ by $\ty$ on $abc$, with new vertex $h$ as center.

\item
$G_1 \to G_2$ by $\ty$ on $ade$, with new vertex $i$ as center.

\item
$G_2 \to G_3$ by $\ty$ on $a\!f\!g$, with new vertex $j$ as center.

\item
$G_3 \to G_4$ by $\ty$ on $bdf$, with new vertex $k$ as center.

\item
$G_4 \to G_5$ by $\ty$ on $beg$, with new vertex $l$ as center.

\item
$G_5 \to G_6$ by $\yt$ on $hi\!j$, deleting vertex $a$.

\item
$G_6 \to G_7$ by $\yt$ on $hkl$, deleting vertex $b$.

}

\section{$G_7$ is not IK}

\thm{
The $\yt$ move does not preserve intrinsic knottedness.
}

\begin{figure}[ht]

 \centering
 \includegraphics[width=80mm]{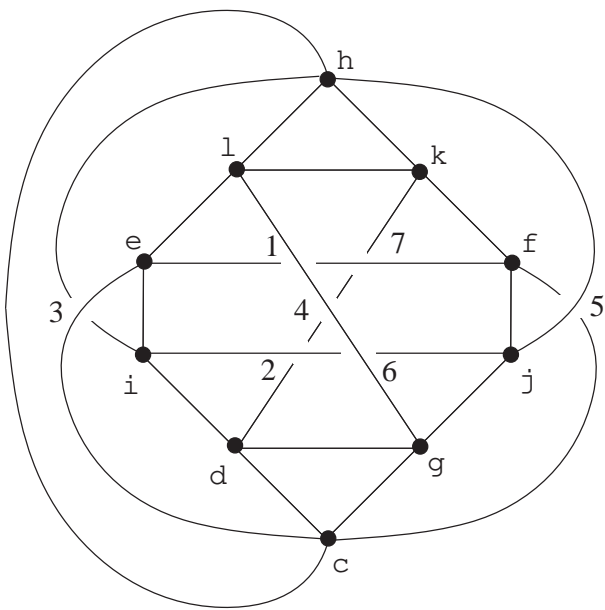}
 \caption{A knotless embedding of $G_7$.
 \label{G7Figure}
 }
\end{figure}

\pf{
Recall that $K_7$ is IK,
$\ty$ moves preserve IKness,
and $G_7$ is obtained from $K_7$
by $\ty$ and $\yt$ moves.
Thus it suffices to prove that
the embedding of $G_7$
shown in Figure~\ref{G7Figure}
has no nontrivial knots.

Figure~\ref{G7Figure} contains seven crossings, numbered 1-7.
Note that rotating this diagram by $180^\circ$
about a horizontal line through its center
leaves the embedded graph invariant,
swaps crossing~1 with 2, and 6 with 7,
and leaves crossings~3, 4, 5 fixed.
And rotating the diagram by $180^\circ$
about a vertical line through the center
also leaves the embedded graph invariant,
but swaps crossing~1 with 7, 2 with 6, and 3 with 5.

Suppose towards contradiction that
this embedded graph contains a nontrivial knot $K$.
The proof consists of the following three steps.

\medskip

\emph{Step 1.}
$K$ must contain exactly one of the edges
$ef$ and $ij$.

Proof:
We will show that if $K$ contains neither or both edges,
then it is a trivial knot.

Suppose $K$ contains neither $ef$ nor $ij$.
Then it does not contain any of the crossings~1, 2, 6, and 7,
and must therefore contain crossings~3, 4, and 5.
Hence $K$ contains the edges
$ec$, $ih$, $dk$, $gl$, $fc$, and $jh$.
If $K$ contains $ei$ or $fj$,
then at least one of its crossings can be untwisted,
making $K$ trivial.
So $K$ must contain $el$, $fk$,
$id$, and $jg$.
Then $K$ is easily seen to be trivial.

Now suppose $K$ contains both $ef$ and $ij$.
Then it cannot contain both 3 and 5,
since otherwise it would be a link.
So, by symmetry,
we can assume $K$ does not contain 5.
Furthermore, if $K$ contains $fj$,
then it is trivial.
It follows that $K$ must contain
at least one of $fk$ or $jg$.
By symmetry,
we can assume it contains $fk$.
We claim that $K$ must contain $dk$,
since otherwise
it will contain at most three crossings, 3, 1, and 6;
but 1 and 6 do not alternate,
which makes $K$ trivial.
Now, $dk\!f\!e$ can be isotoped,
with fixed endpoints,
to eliminate 1, 4, and 7.
So $K$ must contain 3, 2, and 6.
If $K$ contains $jh$,
3 and 6 will not alternate,
making $K$ trivial.
So $K$ must contain $jg$.
But then 6 can be isotoped away,
again making $K$ trivial.
This proves Step~1.

\medskip

So, by symmetry,
we can assume $K$ contains $ef$ and not $ij$.
Hence $K$ does not contain
crossings~2 or 6.

\medskip

\emph{Step 2.}
$K$ must contain 1, 4, and 7.

Proof:
Suppose, towards contradiction,
that $K$ does not contain $gl$.
Then it contains at most three crossings, 3, 5, and 7;
but $ec$, $cf$, and $fe$ form a cycle,
and therefore only links contain all
three crossings 3, 5, and 7.
Hence $K$ contains $gl$.
By a symmetric argument,
$K$ contains $dk$.
Thus $K$ contains crossings~1, 4, and 7.

\medskip

\emph{Step 3.}
$K$ contains exactly one of 3 and 5.

Proof:
If it contains both, it will be a link.
If it contains neither, it will be trivial,
since 1 and 4 do not alternate.

\smallskip

So, by symmetry, we can assume that
$K$ contains 1, 3, 4, and 7,
and no other crossings.
As $K$ does not contain $ij$,
this implies that
$K$ contains $di$.
But $hidk$ is isotopic, with fixed endpoints, to $hk$.
Thus $K$ is isotopic to a knot
that contains only crossing~1,
and therefore is trivial.

}

{\sc Acknowledgments}.
The second author thanks Caltech for its hospitality
while he worked on this paper during his sabbatical leave.

\vspace{3mm}

\vspace{4mm}

{\scriptsize

\noindent
Erica Flapan:
Mathematics Department,
Pomona College,
Claremont, CA 91711, USA.
eflapan@pomona.edu.

\vspace{2mm}

\noindent
Ramin Naimi:
Mathematics Department,
Occidental College,
Los Angeles, CA 90041, USA.
rnaimi-at-oxy.edu.

}

}\end{document}